\documentclass[a4paper,11pt,english,reqno]{amsart}
\usepackage[utf8]{inputenc}
\usepackage[english]{babel}
\usepackage[T1]{fontenc} 
\usepackage{graphicx}
\usepackage[a4paper, margin=2.5cm]{geometry}
\usepackage{amsmath,amssymb,amsfonts,dsfont,amsthm}    
\newcommand{\E}{{\mathcal E}}

\newcommand{\tr}{{\mbox{\rm tr}}}
\newcommand{\mO}{\mathcal{O}}
\newcommand{\C}{\mathds{C}}

\newcommand{\defeq}{\stackrel{\mathrm{def}}{=}}
\newcommand{\prob}{\mathbf{P}}
\newcommand{\erw}{\mathds{E}}
\newtheorem{thm}{Theorem}

\newtheorem{lem}[thm]{Lemma}
\newtheorem{ass}[thm]{Assumption}
\newtheorem{rem}[thm]{Remark}

\numberwithin{equation}{section}
\usepackage{hyperref}
\hypersetup{pdfborder=0 0 0, 
	    colorlinks=true,
	    citecolor=blue,
	    linkcolor=blue,
	    urlcolor=blue,
	    pdfauthor={Martin Vogel, Ofer Zeitouni}
	   }
\title{Deterministic equivalence for noisy perturbations }
\author{Martin Vogel}
\author{Ofer Zeitouni}
\address[Martin Vogel]{Institut de Recherche Math{\'e}matique Avanc{\'e}e - UMR 7501, 
Universit{\'e} de Strasbourg et CNRS, 7 rue René-Descartes, 67084 Strasbourg Cedex, France.}
\address[Ofer Zeitouni]{Department of Mathematics, 
Weizmann Institute of Science, POB 26, Rehovot 76100, Israel
and Courant Institute, New York University,
251 Mercer St, New York, NY 10012, USA.}
 \date{January 23, 2020}
\begin{document}
\maketitle

\begin{abstract}
We prove a quantitative deterministic equivalence theorem for the logarithmic potentials of deterministic complex 
$N\times N$ matrices subject to small random perturbations. We show that with probability close to 
$1$ this log-potential is, up to a small error, determined by the singular values of the unperturbed 
matrix which are larger than some small $N$-dependent cut-off parameter. 
 \end{abstract}
 \setcounter{tocdepth}{1}
%
\section{Introduction and statement of results}
In evaluating the limit of empirical measures of eigenvalues of (non-Hermitian) matrices, an important role is played by the evaluation of certain determinants. Specifically,
for a sequence of matrices $X_N$ of dimension $N\times N$ having eigenvalues $\lambda_i(X_N)$, let
$L_N(X_N)=N^{-1} \sum_{i=1}^N \delta_{\lambda_i(X_N)}$ denote the empirical measure of eigenvalues of $X_N$ and let
$\mathcal{L}_{X_N}(z)=\int \log|z-x|L_N(X_N)(dx)$ denote its log-potential. Since a.e. convergence of log-potentials implies the weak convergence of the associated measures,
the evaluation of limits of log-potentials has played an important role in the study of convergence of the spectrum of random matrices. We refer to \cite{Tao,BC} for introductions to this vast topic.

Since
$$\mathcal{L}_{X_N}(z)= \frac12 \log \det(z-X_N)( z-X_N)^*=\log |\det(z-X_N)| ,$$
 evaluating logarithmic potentials amounts to computing determinants. 
In their study of the spectrum of small, noisy  perturbations of  non-normal matrices, 
the authors of \cite{BPZ} have identified a certain \textit{deterministic equivalent} result, which we now present. 
\begin{thm}\cite[Theorem 2.1]{BPZ}
\label{theo-BPZ}
Let $A=A_N$ be a sequence of
deterministic complex $N\times N$-matrix
of uniformly  bounded norm and
 singular values $s_1\geq \dots s_N \geq 0$. 
Fix $\gamma>1/2$ and $\eta>0$. Set $\varepsilon_N=N^{-\eta}$, and set $N^*$ to be the largest integer $i$ so that
\begin{equation}\label{eq:N*-dfn}
s_{N-i+1}\leq  \varepsilon_N^{-1} N^{-\gamma}(N-i+1)^{1/2}.
\end{equation}
If no such $i$ exists 
then
set  $N^*=1$.  Let $G_N$ be a matrix whose entries are i.i.d. standard complex Gaussian variables. 
 Then, if  $N^* \log N / N \to \alpha < \infty$,
  \begin{equation}\label{eq:log-det-conv}
    \frac{1}{N} \log| \det( A_N + N^{-\gamma}G_N)|
    -
    \frac{1}{N} \sum_{i=1}^{N-N^*+1}  \log s_i    \to_{N\to\infty}  0, \quad \mbox{ in probability}.
  \end{equation}
in probability, as $N\to\infty$. If $\alpha = 0,$ we may take $\varepsilon_N = N^{-\eta}$ for any $\eta > 0.$
\end{thm}

The proof in \cite{BPZ} uses in an essential way the
unitary  invariance of $G_N$, and probabilistic arguments. However, it does not directly extend to other noise models, not even 
to the case where $G_N$ is a matrix consisting of independent real standard Gaussian variables. The purpose of this note is to present a very general version 
of Theorem \ref{theo-BPZ}, based on the Grushin problem studied in 
\cite{Vo19c}. It will be stated under the following assumption on the noise matrix.
Here and throughout, 
for a matrix $A$,
		$s_1(A)\geq s_2(A)\geq \cdots\geq s_N(A)\geq 0$ denote the  singular values of $A$, 
		and
$\|A\|$ denotes the operator norm of $G$, i.e. $\|A\|=s_1(A)$,
\begin{ass}
\label{ass-G}
 $G=G_N$ is  an $N\times N$ random matrix such that the following hold. 
\begin{enumerate}
	\item\textbf{Norm bound} There exists a $\kappa_1 >0 $ such that 
		\begin{equation}\label{rn1.0}
		 	\erw[ \|G\| ] = \mO(N^{\kappa_1}).
		\end{equation}
	\item\textbf{Anti-concentration bound} For each $\theta>0 $ there 
		exists a $\beta >0 $ such that for any fixed deterministic complex 
		$N\times N$ matrix $D$ with  $\|D\| = \mO(N^{\kappa_2})$, $\kappa_2\geq 0$, we have that 
		\begin{equation}\label{rn1.1}
		  \prob (s_N(D + G) \leq N^{-\beta}) = \varepsilon_N (\theta)=o(1).
		\end{equation}
\end{enumerate}
\end{ass}
\begin{thm}
  \label{theo-main}
Let $A=A_N$ be a deterministic complex $N\times N$-matrix with $\|A\|=\mO(N^{\kappa_2})$
for some fixed $\kappa_2\geq 0$, and assume $G=G_N$ satisfies Assumption \ref{ass-G}. Let 
$s_1\geq \dots s_N \geq 0$ denote the singular values of $A$. 
Suppose that for some fixed $L>0 $ there exists  
\begin{equation}
  \label{gp3}
  CN^{-L}\leq \alpha \leq 1
\end{equation}
such that 
\begin{equation}
	 \label{gp4}
	\#\{j; s_j \in [0,\alpha]\} \leq 
	\nu_N \frac{N}{\log N}=:M, \quad \nu_N=o(1).
\end{equation}
For  $\tau>0$ and any fixed $\gamma \gg 1$ we let 
\begin{equation*}
	N^{-\gamma} \leq \delta \ll N^{-\kappa_1}\alpha \tau^{-1}.
\end{equation*}
Then, we have that 
\begin{equation*}
	\left| \frac{1}{N}\log |\det (A + \delta G) | - \frac{1}{N} \sum_{ j: ~s_j >\alpha } \log s_j \right|
	= \mO(1) \left( \nu_N  + \alpha^{-1}N^{\kappa_1} \delta \tau\right)
\end{equation*}
with probability $\geq 1 -\varepsilon_N(\kappa_2+\gamma)- \tau^{-1}$. 
\end{thm}
\begin{rem}
  Assumption \ref{ass-G} holds for a large class of noise matrices, including
  those with iid entries of zero mean and finite variance. We refer to 
  \cite[Remark 1.3]{BPZ1} for details and references.
\end{rem}
\paragraph{\textbf{Acknowledgments}} 
The authors are grateful to the Oberwolfach conference \emph{Random Matrices} in December 2019 
during which most of this text has been written. M. Vogel was supported by the a CNRS 
Momentum 2017 grant. O. Zeitouni was supported by the
European Research Council (ERC) under the European Union's Horizon 2020 research and innovation programme (grant agreement No. 692452). 

\section{Grushin problem}
We now present the proof of Theorem
\ref{theo-main}, based on \cite{Vo19c,HaSj08}, see also \cite{SjVo19a,SjVo19b}
We begin by setting up a well-posed Grushin problem. Let $A=A_N$ 
be a deterministic complex $N\times N$-matrix and let 
 \begin{equation}\label{gp1}
	0\leq t_1^2 \leq \dots \leq t_{N}^2 
\end{equation}
denote the eigenvalues of $A^*A$ with associated orthonormal basis 
of eigenvectors $e_1,\dots,e_{N}\in \C^N$. The spectra of $A^*A$ 
and $AA^*$ are equal and we can find an orthonormal basis 
$f_1,\dots,f_{N}\in \C^N$ of eigenvectors of $AA^*$ associated with the 
eigenvalues \eqref{gp1} such that 
 \begin{equation}\label{gp2}
	A^* f_i = t_i e_i, \quad Ae_i = t_i f_i, \quad i=1,\dots, N.
\end{equation}
Recall $\alpha,M$, see \eqref{gp3},\eqref{gp4}, and let $\delta_i$, $1\leq i \leq M$, 
denote an orthonormal basis of $\C^M$.
Put 
\begin{equation}\label{gp5}
R_+=\sum_{i=1}^M \delta_i\circ  e_i^*, \quad  R_-=\sum_{i=1}^M f_i\circ\delta_i^*,
\end{equation}
We claim that  the Grushin problem 
\begin{equation}\label{gp6}
\mathcal{P} = \begin{pmatrix}
			A & R_- \\ R_+ & 0 \\
		\end{pmatrix} : \C^N\times \C^M \longrightarrow  \C^N\times \C^M 
\end{equation}
is bijective. To see this we take $(v,v_+)\in \C^N\times \C^M$ and we want to solve
 \begin{equation}\label{gp7}
		\mathcal{P}
		 \begin{pmatrix} u \\ u_- 
		\end{pmatrix} =
		 \begin{pmatrix} v \\ v_+
		\end{pmatrix}.
\end{equation}
We write $u= \sum_1^{N} u(j) e_j$ and $v= \sum_1^{N} v(j) f_j$. Similarly, we express 
$u_-,v_+$ in the basis $\delta_1,\dots,\delta_M$. 
The relation \eqref{gp2} then shows 
that \eqref{gp7} is equivalent to 
\begin{equation*}
		\begin{cases}
			\sum_1^{N} t_i u_i f_i + \sum_1^M u_-(j)f_j = \sum_1^{N} v_j f_j \\ 
			u_j = v_+(j), \quad j =1, \dots, M,
		\end{cases}
\end{equation*}
which can be written as 
\begin{equation}\label{gp7.1}
		\begin{cases}
		 \qquad  \quad t_i u_i f_i  = 
		v_i f_i, \qquad \qquad\qquad i=M+1,\ldots,N,  \\[1.5ex] 
			 \begin{pmatrix} t_i & 1 \\ 1 & 0 \\
			\end{pmatrix}
			 \begin{pmatrix} u_i \\ u_-(i)\\
			\end{pmatrix}
			= 
			 \begin{pmatrix} v_i \\ v_+(i)
			\end{pmatrix}, \quad i=1,\dots, M.
		\end{cases}
\end{equation}
Since 
\begin{equation*}
	 \begin{pmatrix} t_i & 1 \\ 1 & 0 \\
	\end{pmatrix}^{-1} 
	=
	 \begin{pmatrix} 0& 1 \\ 1 & -t_i  \\
	\end{pmatrix},
\end{equation*}
we see that 
\begin{equation}\label{gp7.2}
		\mathcal{P}^{-1}= \mathcal{E} =
		\begin{pmatrix} 
		E  & E_+\\ 
		E_- & E_{-+}\\
		\end{pmatrix} 
\end{equation}
with
\begin{equation}\label{gp8}
\begin{split}
		&E = \sum_{M+1}^{N} \frac{1}{t_i} e_i \circ f_i, \quad 
	         E_+ = \sum_1^M e_i \circ \delta_i^*,  \\
		&E_- =\sum_1^M  \delta_i\circ f_i^*, \quad
		 E_{-+} = - \sum_1^M t_j \delta_j\circ\delta_j^*, 
\end{split}
\end{equation}
and the norm estimates
\begin{equation}\label{gp9}
	\|E(z) \| \leq \frac{1}{\alpha}, \quad \| E_{\pm } \| =1, \quad 
	\| E_{-+}\| \leq \alpha.
\end{equation}
Furthermore, \eqref{gp7.1} shows that 
\begin{equation}\label{gp10}
	|\det\mathcal{P}|^2 =\prod_{M+1}^N t_i^2.
\end{equation}
\subsection{Gruhsin problem for the perturbed operator}
Now we turn to the perturbed operator 
\begin{equation}\label{gpp1}
	A^\delta=A+\delta G, \quad 0 \leq \delta \ll 1.
\end{equation}
where $G$ is a complex $N\times N$-matrix. Let $R_{\pm}$ be 
as in \eqref{gp5}, and put 
\begin{equation}\label{gpp2}
\mathcal{P}^{\delta} = \begin{pmatrix}
			A^{\delta} & R_- \\ R_+ & 0 \\
		\end{pmatrix} : \C^N\times \C^M \longrightarrow  \C^N\times \C^M 
\end{equation}
Then $\mathcal{P}=\mathcal{P}^0$. Applying $\mathcal{E}$, 
see \eqref{gp7.2},  from
the right to \eqref{gpp2} yields 
\begin{equation}\label{gpp3}
  \mathcal{P}^{\delta}\mathcal{E} = I_{N+M} + 
\begin{pmatrix}
			\delta G E &  \delta G E_+\\ 0 & 0 \\
		\end{pmatrix} 
\end{equation}
Suppose that 
\begin{equation}\label{gpp4}
	\delta \|G\| \alpha^{-1} \leq \frac{1}{2}.
\end{equation}
Then, see \eqref{gpp1}, the matrix
$\mathcal{P}^{\delta}\mathcal{E} $ is invertible by a Neumann series 
argument and we get that 
\begin{equation}\label{gpp5}
\begin{split}
	\mathcal{E}^{\delta} = (\mathcal{P}^{\delta})^{-1} 
	&= \mathcal{E} + \sum_{n=1}^{\infty}(-\delta)^n 
	\begin{pmatrix}
			E(G E)^n &   (EG)^{n}E_+\\  E_-(G E)^n& E_-(GE)^{n-1}GE_+ \\
	\end{pmatrix} \\
	& \defeq \begin{pmatrix} 
		E^{\delta}  & E^{\delta}_+\\ 
		E^{\delta}_-& E^{\delta}_{-+}\\
		\end{pmatrix},
\end{split}
\end{equation}
where by \eqref{gpp4}, \eqref{gp9}, 
 \begin{equation}\label{gpp6}
 \begin{split}
		&\| E^{\delta} \| = \|  E( 1+ \delta GE)^{-1} \| \leq 2 \|E\| \leq 2 \alpha^{-1}, \\
		&\| E_+^{\delta} \| = \|  ( 1+ \delta GE)^{-1}E_+ \| \leq 2 \|E_+\| \leq 2,  \\
		&\| E_-^{\delta} \| = \| E_- ( 1+ \delta GE)^{-1} \| \leq 2\|E_-\| \leq 2,  \\
		&\| E_{-+}^{\delta} -E_{-+}\| = 
		\| E_- ( 1+ \delta Q_{\omega}E)^{-1}\delta GE_+ \| \leq 2 \|\delta G \| \leq  \alpha.  \\
\end{split}
\end{equation}
The Schur complement formula applied to $\mathcal{P}^{\delta}$ and 
$\mathcal{E}^{\delta}$ shows that 
 \begin{equation}\label{gpp7}
	\log |\det A^\delta|=\log |\det \mathcal{P}^\delta|+\log |\det E_{-+}^\delta|.
\end{equation}
Notice that 
 \begin{equation}\label{gpp8}
	\begin{split}
		\left|  \log |\det \mathcal{P}^\delta|-\log |\det \mathcal{P}^0| \right| &
		  =	\left|\Re \int_0^\delta  \tr (\E^\tau
		\frac{d}{d\tau} \mathcal{P}^\tau) d\tau\right|\\
		&=	\left|\Re \int_0^\delta \tr \begin{pmatrix} E^\tau&E_+^\tau\\
  		E_-^\tau& E_{-+}^\tau\end{pmatrix}
		\cdot \begin{pmatrix} G &0\\0&0\end{pmatrix}
		d\tau\right|\\
		&=	\left| \Re \int_0^\delta  \tr(E^\tau G)d\tau \right|\\
		&\leq 2 \alpha^{-1}\delta N \|G\|.
	\end{split}
\end{equation}
Here in the last line we used \eqref{gpp6}. Thus,
\begin{equation}\label{gpp9}
  \left|\frac{1}{N} \log |\det \mathcal{P}^\delta|-\frac{1}{N} \log |\det  \mathcal{P}|\right|
  \leq 2 \alpha^{-1} \delta \|G\|.
\end{equation}
Notice that by \eqref{gpp6}, \eqref{gp9}, we have that  $\|E_{-+}^{\delta}\| \leq 2\alpha$. 
Thus, by \eqref{gpp7} and \eqref{gpp9},
\begin{equation}\label{gpp10}
\log |\det A^\delta| \leq \log |\det \mathcal{P}|+ M |\log 2\alpha|  + 2\alpha^{-1}\delta N \|G\|.
\end{equation}
\subsection{Random noise matrix} 
We recall Assumption \ref{ass-G} on the noise matrix. By Markov's inequality, 
\begin{equation}\label{rn2}
  \prob ( \|G\| >  C N^{\kappa_1} \tau ) \leq \tau^{-1}, \quad \tau >0. 
\end{equation}
Since
\begin{equation}\label{rn3}
	0 \leq \delta \ll N^{-\kappa_1}\alpha \tau^{-1},
\end{equation}
we obtain that 
with probability $\geq 1 - \tau^{-1}$, we have that \eqref{gpp4} holds. Hence, the estimates 
\eqref{gpp6} and \eqref{gpp7} hold with the same probability. This together with 
\eqref{gpp10}, \eqref{gp4} and \eqref{gp3}, implies that $\|G\| \leq CN^{\kappa_1}\tau$ and 
\begin{equation}\label{rn4}
\log |\det A^\delta| \leq \log |\det \mathcal{P}|+ \mO(1) \nu_N N 
	+ \mO(1)\alpha^{-1}N^{1+\kappa_1} \delta \tau
\end{equation}
with probability $\geq 1 - \tau^{-1}$. 
\par
It remains to find a lower bound on $\log |\det E_{-+}^{\delta}|$. We begin by 
recalling a classical result on Grushin problems, see for instance \cite[Lemma 18]{Vo19c}.
\begin{lem}\label{gpp:lem1}
	Let $\mathcal{H}$ be an $N$-dimensional complex Hilbert spaces, and 
	let $N\geq M>0$. Suppose that 	
	 \begin{equation*}
		\mathcal{P} = \begin{pmatrix} 
		P & R_- \\ 
		R_+ & 0\\
		\end{pmatrix} :  \mathcal{H} \times \C^M \longrightarrow  \mathcal{H} \times \C^M
	\end{equation*}
	is a bijective matrix of linear operators, with inverse 
	 \begin{equation*}
		\mathcal{E}
		=
		 \begin{pmatrix} 
		E  & E_+\\ 
		E_-& E_{-+}\\
		\end{pmatrix}.
	\end{equation*}
	Let $0 \leq t_1(P) \leq \dots \leq t_N(P)$ denote the eigenvalues of $(P^*P)^{1/2}$, and 
	let $0 \leq t_1(E_{-+}) \leq \dots \leq t_M(E_{-+})$ denote the eigenvalues of $(E_{-+}^*E_{-+})^{1/2}$. 
	Then, 
	\begin{equation*}
			\frac{t_n(E_{-+})}{\| E\| t_n(E_{-+}) + \|E_-\| \|E_+\|}
			\leq  t_n(P) \leq \|R_+\| \|R_-\| t_n(E_{-+}), \quad 1\leq n \leq M.
	\end{equation*}
\end{lem}
By \eqref{gp5} we know that $\|R_{\pm}\| =1$, and by \eqref{gpp6} we then 
get 
\begin{equation}\label{rn6}
		 t_n(A^{\delta}) \leq  t_n(E^{\delta}_{-+}), \quad 1\leq n \leq M.
\end{equation}
Next note that, for any $\delta\geq N^{-\gamma}$ and $\beta>0$
and any deterministic matrix $A$,
  \begin{eqnarray*}
    \textbf{P}\left( s_{N}( A + \delta G) \leq  N^{-\gamma -\beta} \right) 
    &=&
    \textbf{P}\left( s_{N}( A/\delta +  G) \leq  N^{-\gamma -\beta}/\delta
    \right)\\
    &\leq &
    \textbf{P}\left( s_{N}( A/\delta +  G) \leq  N^{-\beta}
    \right).
  \end{eqnarray*}
Thus, from \eqref{rn1.1}, 
there exists a $\beta >0$ such for any fixed deterministic matrix $A$ with 
$\|A \| = \mO(N^{\kappa_2})$ and any $\delta \geq N^{-\gamma}$,
we have that 
	\begin{equation}\label{r7}
	\textbf{P}\left( s_{N}( A + \delta G) \leq  N^{-\gamma -\beta} \right) 
	\leq 
      \varepsilon_N(\kappa_2+\gamma).
	\end{equation}	
We recall that
\begin{equation}\label{rn5}
	N^{-\gamma} \leq \delta \ll N^{-\kappa_1}\alpha \tau^{-1}.
\end{equation}
By combining \eqref{rn6}, \eqref{r7} and \eqref{rn2}, we obtain that 
\begin{equation}\label{rn7}
	\textbf{P}\left( s_{M}( E_{-+}^{\delta}) >  N^{-\gamma-\beta} 
	 \text{ and } \Vert G \Vert \leq CN^{\kappa_1} \tau \right) 
	 \geq 1 - \varepsilon_N(\kappa_2+\gamma) - \tau^{-1}.
\end{equation}	
Provided that this event holds and using \eqref{gp4}, we get hat 
\begin{equation}\label{rn8}
\begin{split}
\log |\det E_{-+}^{\delta} | & = \sum_1^M \log s_j(E_{-+}^{\delta}) \\ 
& \geq M \log s_M(E_{-+}^{\delta})  \\
& \geq - (\gamma + \beta) M \log N \\
&\geq - (\gamma + \beta) \nu_N N,
\end{split}
\end{equation}
which in combination with \eqref{gpp7}, \eqref{gpp9} yields that 
 \begin{equation}\label{rn9}
	\log |\det A^\delta| 
	\geq 
	\log |\det \mathcal{P}| -  \mO(1) \nu_NN 
	- \mO(1)\alpha^{-1}N^{1+\kappa_1} \delta \tau
\end{equation}
with probability $\geq 1 - \varepsilon_N(\kappa_2+\gamma) - \tau^{-1}$. This, in view of \eqref{rn4}, \eqref{rn3} 
and \eqref{gp10} concludes the proof of the theorem. 
\providecommand{\bysame}{\leavevmode\hbox to3em{\hrulefill}\thinspace}
\providecommand{\MR}{\relax\ifhmode\unskip\space\fi MR }
\providecommand{\MRhref}[2]{%
  \href{http://www.ams.org/mathscinet-getitem?mr=#1}{#2}
}
\providecommand{\href}[2]{#2}

\end{document}